# LOWEST TERMS REVISITED

Laurence D. Robinson
Department of Mathematics
Ohio Northern University
Ada, OH 45810
L-Robinson.1@onu.edu

Carter G. Lyons
Department of Mathematics and Statistics
James Madison University
Harrisonburg, VA 22807

**Abstract**: In the September 1994 issue of Math Horizons the following problem is given in the "Problem Section" (p. 33): "Problem 5: Lowest Terms – What fraction has the smallest denominator in the interval $\left(\frac{19}{94}, \frac{17}{76}\right)$?" In this paper we develop a general algorithm which gives a systematic procedure for solving problems of this type.

We consider the interval $\left(\frac{A_0}{B_0}, \frac{C_0}{D_0}\right)$, in which $A_0/B_0$ and $C_0/D_0$ are reduced proper fractions. We wish to find the reduced proper fraction $P_0/Q_0$ falling inside this interval which has the minimum possible denominator $Q_0$.

We note that there exists a unique value of $P_0$ associated with the minimum $Q_0$ value, and that this $P_0$ value is the minimum possible value for the numerator of any reduced proper fraction within the interval. This can be seen as follows: Consider the two reduced proper fractions $P_0/Q_0$ and $P^*/Q^*$, with $P^* < P_0$ and $Q^* \geq Q_0$. It readily follows that

$$\frac{P^*}{Q^*} < \frac{P^*}{Q_0 - (P_0 - P^*)} < \frac{P_0}{Q_0}.$$

Thus, there exists a reduced proper fraction falling between $P^*/Q^*$ and $P_0/Q_0$ having a denominator which is smaller than $Q_0$. Now suppose that $P_0/Q_0$ is the reduced proper fraction falling inside the interval $\left(\frac{A_0}{B_0}, \frac{C_0}{D_0}\right)$ having the minimum possible denominator. It follows that there cannot exist a reduced proper fraction $P^*/Q^*$, with $P^* < P_0$ and $Q^* \geq Q_0$, which falls inside the interval. We thus conclude that the numerator $P_0$ associated with the minimum possible denominator $Q_0$ is unique, and is the minimum possible numerator for the interval.

To determine $P_0/Q_0$ we use an iterative procedure which works as follows: We express the original interval $\left(\frac{A_0}{B_0}, \frac{C_0}{D_0}\right)$ as $\left(\frac{1}{W_1 + C_1/D_1}, \frac{1}{W_1 + A_1/B_1}\right)$, in such a way that $A_1/B_1$ is either a reduced proper fraction or 0. We then form the new interval $\left(\frac{A_1}{B_1}, \frac{C_1}{D_1}\right)$, which in turn is expressed as $\left( \dfrac{1}{W_2 + C_2/D_2} \;,\; \dfrac{1}{W_2 + A_2/B_2} \right)$, in such a way that $A_2/B_2$ is either a reduced proper fraction or 0, and so forth. The process is repeated until a terminal interval $\left( \dfrac{A_k}{B_k}, \dfrac{C_k}{D_k} \right)$ is obtained such that $A_k/B_k$ is either a reduced proper fraction or 0, and $C_k/D_k > 1$. This terminal form is guaranteed to exist by the division algorithm. For such an interval it is clear that the reduced fraction $P_k/Q_k = 1$ has minimum possible values for both the denominator and numerator, specifically
$P_k = 1$ and $Q_k = 1$.

Now letting $P_k = 1$ and $Q_k = 1$, we determine a reduced proper fraction $P_{k-1}/Q_{k-1} = \dfrac{1}{W_k + P_k/Q_k} = \dfrac{Q_k}{W_k Q_k + P_k} = \dfrac{1}{W_k + 1}$ . Since $P_k/Q_k$ falls inside the interval $\left( \dfrac{A_k}{B_k}, \dfrac{C_k}{D_k} \right)$, it follows that $P_{k-1}/Q_{k-1}$ falls inside the interval $\left( \dfrac{A_{k-1}}{B_{k-1}}, \dfrac{C_{k-1}}{D_{k-1}} \right) = \left( \dfrac{1}{W_k + C_k/D_k} \;,\; \dfrac{1}{W_k + A_k/B_k} \right)$. Furthermore, since the values $P_k = 1$ and $Q_k = 1$ are both minimal, we see that $P_{k-1}/Q_{k-1}$ is the reduced proper fraction having both minimal numerator and denominator inside the interval $\left( \dfrac{A_{k-1}}{B_{k-1}}, \dfrac{C_{k-1}}{D_{k-1}} \right)$.

In a similar manner we obtain

$P_{k-2}/Q_{k-2} = \dfrac{1}{W_{k-1} + P_{k-1}/Q_{k-1}} = \dfrac{Q_{k-1}}{W_{k-1} Q_{k-1} + P_{k-1}}$, which using the same logic as above is seen to be the reduced proper fraction falling inside the interval $\left( \dfrac{A_{k-2}}{B_{k-2}}, \dfrac{C_{k-2}}{D_{k-2}} \right)$ having both minimal numerator and denominator. The process can be continued until the reduced proper fraction $P_0/Q_0$, having minimal denominator (and numerator) inside the initial interval $\left( \dfrac{A_0}{B_0}, \dfrac{C_0}{D_0} \right)$, is obtained.

In the course of applying the algorithm it is possible to obtain an interval $\left(\frac{A_{k-1}}{B_{k-1}}, \frac{C_{k-1}}{D_{k-1}}\right)$ with $A_{k-1}/B_{k-1} = 0$ and $C_{k-1}/D_{k-1} \le 1$. In such cases proceeding with the algorithm would require division by 0. To avoid this, the procedure can be modified by replacing $A_{k-1}/B_{k-1} = 0$ by $1/(D_{k-1} + 2)$, and then iterating one additional time to obtain the terminal interval $\left(\frac{A_k}{B_k}, \frac{C_k}{D_k}\right)$. The rationale underlying this modification is that the reduced proper fraction $P_{k-1}/Q_{k-1}$ falling inside the interval $\left(\frac{1}{D_{k-1}+2}, \frac{C_{k-1}}{D_{k-1}}\right)$, which has minimal numerator and denominator, must be the same as the reduced proper fraction $P_{k-1}/Q_{k-1}$ falling inside the interval $\left(0, \frac{C_{k-1}}{D_{k-1}}\right)$ having minimal numerator and denominator. This follows from the fact that the desired $P_{k-1}/Q_{k-1}$ must exceed $1/(D_{k-1} + 2)$, since $1/(D_{k-1} + 1)$ falls inside the interval $\left(\frac{1}{D_{k-1}+2}, \frac{C_{k-1}}{D_{k-1}}\right)$ and has a denominator smaller than $1/(D_{k-1} + 2)$.

As an example for which the above modification is not required, consider the initial interval $\left(\frac{A_0}{B_0}, \frac{C_0}{D_0}\right) = \left(\frac{19}{94}, \frac{17}{76}\right)$. To find the reduced proper fraction $P_0/Q_0$ falling inside the interval having minimum denominator $Q_0$ we express the interval as

$$\left(\frac{19}{94}, \frac{17}{76}\right) = \left(\frac{1}{94/19}, \frac{1}{76/17}\right) = \left(\frac{1}{4+\ 18/19}\ ,\ \frac{1}{4+\ 8/17}\right).$$

Next we form the interval $\left(\frac{A_1}{B_1}, \frac{C_1}{D_1}\right) = \left(\frac{8}{17}, \frac{18}{19}\right)$ and repeat the process:

$$\left(\frac{8}{17}, \frac{18}{19}\right) = \left(\frac{1}{17/8}, \frac{1}{19/18}\right) = \left(\frac{1}{1+\ 9/8}\ ,\ \frac{1}{1+\ 1/18}\right).$$

We now form the terminal interval $\left(\frac{A_2}{B_2}, \frac{C_2}{D_2}\right) = \left(\frac{1}{18}, \frac{9}{8}\right)$, which has $C_2/D_2 = 9/8 > 1$. Thus, we obtain $P_2/Q_2 = 1$, from which we obtain $P_1/Q_1 = \frac{1}{1+1} = 1/2$, from which we obtain the desired solution $P_0/Q_0 = \frac{1}{4+1/2} = 2/9 = .2222$. This is the reduced proper fraction $P_0/Q_0$ falling inside the initial interval $\left(\frac{19}{94}, \frac{17}{76}\right) = (.2021, .2237)$ having minimum denominator $Q_0 = 9$.

As an example for which the algorithm requires the described modification, consider the initial interval $\left(\frac{A_0}{B_0}, \frac{C_0}{D_0}\right) = \left(\frac{3}{4}, \frac{17}{22}\right)$. To find the reduced proper fraction $P_0/Q_0$ falling inside the interval having minimum denominator $Q_0$ we express the interval as

$$\left(\frac{3}{4}, \frac{17}{22}\right) = \left(\frac{1}{4/3}, \frac{1}{22/17}\right) = \left(\frac{1}{1+1/3}, \frac{1}{1+5/17}\right).$$

Next we form the interval $\left(\frac{A_1}{B_1}, \frac{C_1}{D_1}\right) = \left(\frac{5}{17}, \frac{1}{3}\right)$ and repeat the process:

$$\left(\frac{5}{17}, \frac{1}{3}\right) = \left(\frac{1}{17/5}, \frac{1}{3}\right) = \left(\frac{1}{3+2/5}, \frac{1}{3+0}\right).$$

Next we form the interval $\left(\frac{A_2}{B_2}, \frac{C_2}{D_2}\right) = \left(0, \frac{2}{5}\right)$, and now we apply the modification, replacing 0 by $1/(D_2+2) = 1/(5+2) = 1/7$. We now repeat the process once more:

$$\left(\frac{1}{7}, \frac{2}{5}\right) = \left(\frac{1}{7}, \frac{1}{5/2}\right) = \left(\frac{1}{2+5}, \frac{1}{2+1/2}\right).$$

We now form the terminal interval $\left(\frac{A_3}{B_3}, \frac{C_3}{D_3}\right) = \left(\frac{1}{2}, 5\right)$, which has $C_3/D_3 = 5 > 1$. Thus, we obtain $P_3/Q_3 = 1$, from which we obtain $P_2/Q_2 = \frac{1}{2+1} = 1/3$, from which we obtain $P_1/Q_1 = \frac{1}{3+1/3} = 3/10$, from which we obtain the desired solution $P_0/Q_0 = \frac{1}{1+3/10} = 10/13 = .7692$.

This is the reduced proper fraction $P_0/Q_0$ falling inside the initial interval $\left(\frac{3}{4}, \frac{17}{22}\right) = (.7500, .7727)$ having minimum denominator $Q_0 = 13$.